\documentstyle[11pt]{article}	%
 \def\be{\begin{equation}}
 \def\ee{\end{equation}}
 
 \def\set{\setcounter{equation}{0}}
 
 \newcounter{equationsave}
 \newenvironment{eqnarrayabc}%
 {\bgroup
 \addtocounter{equation}{1}
 \setcounter{equationsave}{\value{equation}}
 \setcounter{equation}{0}
 \begin{eqnarray}}%
 {\end{eqnarray}
 \setcounter{equation}{\value{equationsave}}\egroup}
 \newcounter{equationsave1}
 \newenvironment{eqnarrayabc*}%
 {\bgroup
 \setcounter{equationsave1}{\value{equation}}
 \setcounter{equation}{1}
 \begin{eqnarray}}%
 {\end{eqnarray}
 \setcounter{equation}{\value{equationsave1}}\egroup}

\begin{document}

\title{On the stability of the Bareiss and related\\
	Toeplitz factorization algorithms%
\thanks{Copyright \copyright\ 1993--2010, the authors.
	\hspace*{\fill} rpb144tr~typeset~using~\LaTeX%
}}

\author{A. W. Bojanczyk\\
School of Electrical Engineering\\
Cornell University\\
Ithaca, NY 14853-5401, USA
	\and
R. P. Brent\\
Computer Sciences Laboratory\\
Australian National University\\
Canberra, ACT 0200, Australia
	\and
F. R. de Hoog\\
Division of Maths.\ and Stats.\\
CSIRO\\
Canberra, ACT 2601, Australia
	\and
D. R. Sweet\\
Electronics Research Laboratory\\
DSTO\\ %
Salisbury, SA 5108, Australia\\[15pt]}

\date{Report TR-CS-93-14\\
	November 1993}

\maketitle
\vspace{-15pt}

\begin{abstract}
This paper contains a numerical stability
analysis of factorization algorithms for
computing the Cholesky decomposition of symmetric positive definite
matrices of displacement rank~2. The algorithms in the class
can be expressed as
sequences of {\it elementary downdating} steps.
The stability of the factorization
algorithms follows directly from the numerical
properties of algorithms for realizing
elementary downdating operations.
It is shown that the Bareiss algorithm for
factorizing a symmetric positive definite Toeplitz matrix
is in the class and hence the Bareiss algorithm is stable.
Some numerical experiments that compare
behavior of the Bareiss algorithm and the Levinson algorithm are presented.
These experiments indicate that in general
(when the reflection coefficients are not all positive)
the Levinson algorithm
can give much larger residuals than the Bareiss algorithm.
\end{abstract}

\thispagestyle{empty}

\section{Introduction}
\label{Sec:1}

We consider the numerical stability of algorithms
for solving a linear system
\be
Tx = b,		\label{eq11}
\ee
where $T$ is an $ n \times n $ positive
definite Toeplitz matrix and $b$ is an $n \times 1 $ vector.
We assume that the system is solved in
floating point arithmetic with relative precision $\epsilon$ by
first computing the Cholesky factor of $T$. Hence the emphasis of the
paper is on factorization algorithms for the matrix $T$.

Roundoff error analyses of Toeplitz systems solvers have been given by
Cybenko~\cite{Cybenko1} and Sweet~\cite{Sweet}. Cybenko showed
that the Levinson-Durbin algorithm produces a residual which, under the
condition that all reflection coefficients are positive, is of comparable
size to that produced by the well behaved Cholesky method.
He hypothesised    %
that the same is true even if the reflection coefficients are not all
positive.
If correct, this would indicate that numerical quality
of the Levinson-Durbin algorithm
is comparable to that of the Cholesky method.

In his PhD thesis~\cite{Sweet}, Sweet presented a roundoff error analysis
of a variant of the Bareiss algorithm~\cite{Bareiss}, and concluded that the
algorithm is numerically stable (in the sense specified in Section~\ref{Sec:7}).
In this paper we strengthen and generalize these early results
on the stability of the Bareiss algorithm. In particular, our approach via
elementary downdating greatly simplifies roundoff error analysis and makes it
applicable to a larger-than-Toeplitz class of matrices.

After introducing the notation and the concept
of {\it elementary downdating} in Sections~\ref{Sec:2} and~\ref{Sec:3},
in Section~\ref{Sec:4} we derive
matrix factorization algorithms as a sequence of elementary
downdating operations (see also~\cite{BBH}).
In Section~\ref{Sec:5} we present a
first order analysis by bounding the first term in an asymptotic expansion
for the error in powers of $\epsilon$.	%
By analyzing the propagation of first order error in the sequence
of downdatings that define the algorithms, we obtain bounds on the
perturbations of the factors in the decompositions.
We show that the computed upper triangular factor $\tilde U$ of a positive
definite Toeplitz matrix $T$ satisfies
\[
T=\tilde U^T\tilde U +\Delta T\;,\;\;\; ||\Delta T|| \leq c(n) \epsilon
||T||\;,	%
\]
where $c(n)$ is a low order polynomial in $n$ and is independent of the
condition number of $T$.
Many of the results of Sections~\ref{Sec:2}--\ref{Sec:5} %
were first reported in~\cite{BBH2},
which also contains some results on the stability of
Levinson's algorithm.

In Section~\ref{Sec:6} we discuss the connection with the Bareiss algorithm and
conclude that the Bareiss algorithm is stable for the class of symmetric
positive definite matrices. Finally, in Section~\ref{Sec:7} we report some
interesting numerical examples that contrast the behaviour of the
Bareiss algorithm with that of the Levinson algorithm.
We show numerically that, in cases where the reflection coefficients
are not all positive, the Levinson algorithm can give much larger residuals
than the Bareiss or Cholesky algorithms.

\section {Notation}
\label{Sec:2}
\set

Unless it is clear from the context, all vectors are real and of
dimension $n$.  Likewise, all matrices are real and their default
dimension is $n \times n$.  If ${\bf a} \in \Re^n$, $\| {\bf a} \|$
denotes the usual Euclidean norm, and if $T \in \Re^{n \times n}$,
$\| T \|$ denotes the induced matrix norm:
\[
\| T \| = \max_{\|{\bf a}\| = 1} \| T {\bf a} \|\;.
\]

Our primary interest is in a symmetric positive definite Toeplitz matrix
$T$ whose $i,j$th entry is
\[
t_{ij} = t_{|i-j|}.
\]

We denote by ${\bf e}_k$, $k=1, \ldots ,n$, the unit vector whose
$k$th element is~$1$ and whose other elements are~$0$.
We use the following special matrices:
\[
Z \equiv \sum_{i=1}^{n-1}e_{i+1} e_i^T =
{
\left( \begin{array}{ccccc} 0 & & \cdots & \cdots & 0 \\
1 & 0 & \cdots & 0 & 0\\
0 & \ddots & & \vdots & \vdots \\
\vdots & & \ddots & 0 & 0 \\
0 & \cdots & 0 & 1 & 0 \\
\end{array} \right) } \;,\\
\]
\[
J\equiv \sum_{i=1}^{n}e_{n-i+1}e_i^T =
\left( \begin{array}{ccccc}0&\cdots & \cdots & 0 & 1\\
\vdots & & \cdot & 1 & 0\\
\vdots & \cdot & \cdot & \cdot & \vdots \\
0 & 1 & \cdot & & \vdots \\
1 & 0 & \cdots & \cdots & 0
\end{array} \right) \;.
\]
The matrix $Z$ is known as a {\it shift-down} matrix.
We also make use of powers of the matrix $Z$, for which
we introduce the following notation:
\[	            		%
Z_k = 	\left\{ \begin{array}{ll}
	I 					& \mbox{if $k=0$,}\\
	Z^k					& \mbox{if $k > 0$.}
	\end{array}
	\right.
\]
The antidiagonal matrix $J$ is called a {\it reversal} matrix,
because the effect of
applying $J$ to a vector is to reverse the order of components of the
vector:
$$
J\left[ \begin{array}{c}x_1\\x_2\\ \vdots \\x_n\end{array} \right]=
\left[ \begin{array}{c}x_n\\x_{n-1}\\ \vdots \\x_1\end{array} \right]\;.
$$

The {\it hyperbolic rotation}
matrix $H(\theta) \in \Re^{2 \times 2}$ is defined by
\be
H(\theta)= \frac{1}{\cos \theta} \left[ \begin{array}{cc}
1 & -\sin \theta\\
- \sin \theta & 1
\end{array} \right]\;.
\label{hyprotH}
\ee
The matrix $H(\theta)$ satisfies the relation
\[
H(\theta)\left[ \begin{array}{cc}
1&0\\0&-1\end{array} \right]H(\theta)=
\left[ \begin{array}{cc}
1&0\\0&-1\end{array} \right]\;,
\]
and it has eigenvalues $\lambda_1(\theta)$, $\lambda_2(\theta)$ given by
\be
\lambda_1(\theta) = \lambda^{-1}_2(\theta) = \sec \theta - \tan \theta.
\label{Heigenval}
\ee
For a given pair of real numbers $a$ and $b$ with $|a|>|b|$, there exists
a hyperbolic rotation matrix $H(\theta)$ such that
\be
H(\theta) \left[ \begin{array}{c}a\\b\end{array} \right]=
\left[ \begin{array}{c}\sqrt{a^2-b^2} \\0\end{array} \right]\;.
\label{hyprotzero}
\ee
The angle of rotation $\theta $ is determined by
\be
\sin \theta = \frac{b}{a} \;\;\;,\;\;\; \cos \theta =
\frac{\sqrt{a^2-b^2}}{a}\;\;.
\label{hyprotangle}
\ee

\section {Elementary Downdating}
\label{Sec:3}
\set

In this section we introduce the concept of elementary downdating.
The elementary downdating problem is a special case of
a more general downdating problem that arises
in Cholesky factorization of a positive definite difference of
two outer product matrices~\cite{Alexander,BBDH,Bojanczyk,Fletcher}.
In Section~\ref{Sec:4}, factorization algorithms are derived
in terms of a sequence of downdating steps. The numerical properties
of the algorithms are then related to the properties of the sequence of
elementary downdating steps.

\pagebreak[3]
Let ${\bf u}_k$, ${\bf v}_k \in \Re^n$ have the following form:
{%
\begin{center}
\begin{tabular}{cccccccccc}
& & & & & $k$ & & & & \\
& & & & & $\downarrow$ & & & & \\
${\bf u}_k^T$ & = & [0 & $\ldots$ & 0 & $\times$ & $\times$ & $\times$
& $\ldots$ & $\times$ ]$\;,$ \\
${\bf v}_k^T$ & =&
[0&$\ldots$&0&0&$\times$&$\times$&$\ldots$&$\times$ ]$\;,$ \\
&&&&&&$\uparrow$&&&\\
&&&&&&$k+1$&&&\\
\end{tabular}
\end{center}
}
\noindent
that is:
\[
{\bf e}_j^T{\bf u}_k=0 \; , \; j<k\;,\;\;\;{\rm and}\;\;\;\;
{\bf e}_j^T{\bf v}_k=0 \; , \; j \leq k\;.
\]
Applying the shift-down matrix $Z$ to ${\bf u}_k$, we have
{%
\begin{center}
\begin{tabular}{cccccccccc}
&&&&&&$k+1$&&&\\
&&&&&&$\downarrow$&&&\\
${\bf u}_k^TZ^T$ & =&
[0&$\ldots$&0&0&$\times$&$\times$&$\ldots$&$\times$] ,\\
${\bf v}_k^T$ & =&
[0&$\ldots$&0&0&$\times$&$\times$&$\ldots$&$\times$] .\\
&&&&&&$\uparrow$&&&\\
&&&&&&$k+1$&&&\\
\end{tabular}
\end{center}
}
Suppose that we wish to find ${\bf u}_{k+1}$, ${\bf v}_{k+1} \in \Re^n$
to satisfy
\be
{\bf u}_{k+1}{\bf u}_{k+1}^T - {\bf v}_{k+1}{\bf v}_{k+1}^T =
Z{\bf u}_k{\bf u}^T_k Z^T-{\bf v}_k{\bf v}_k^T,
\label{elementdown}
\ee
where
{%
\begin{center}
\begin{tabular}{cccccccccc}
&&&&&&$k+1$&&&\\
&&&&&&$\downarrow$&&&\\
${\bf u}_{k+1}^T$ & =&
[0&$\ldots$&0&0&$\times$&$\times$&$\ldots$&$\times$] ,\\
${\bf v}_{k+1}^T$ & =&
[0&$\ldots$&0&0&0&$\times$&$\ldots$&$\times$] ,\\
&&&&&&&$\uparrow$&&\\
&&&&&&&$k+2$&&\\
\end{tabular}
\end{center}
}
\noindent
that is
\be
{\bf e}_j^T{\bf u}_{k+1}=0 \; , \; j<k+1\;, \;\;\;\;{\rm and}\;\;\;\;
{\bf e}_j^T{\bf v}_{k+1}=0 \; , \; j \leq k+1\;.
\label{zerok+1}
\ee
We refer to the problem of finding ${\bf u}_{k+1}$ and ${\bf v}_{k+1}$
to satisfy~(\ref{elementdown}), given
${\bf u}_{k}$ and ${\bf v}_{k}$, as the {\it elementary downdating} problem.
It can be rewritten as follows:
\[
\left[{\bf u}_{k+1} \; {\bf v}_{k+1}\right] \left[ \begin{array}{cc}
1 & 0 \\ 0 & -1 \end{array} \right] \left[ \begin{array}{c}{\bf
u}_{k+1}^T\\ {\bf v}_{k+1}^T\end{array} \right] =
\left[Z{\bf u}_{k} \; {\bf v}_{k}\right] \left[ \begin{array}{cc}
1 & 0 \\ 0 & -1 \end{array} \right] \left[ \begin{array}{c}
{\bf u}_{k}^TZ^T\\ {\bf v}_{k}^T\end{array} \right]\;.
\]
From~(\ref{hyprotH}),~(\ref{hyprotzero}) and~(\ref{hyprotangle}),
it is clear that the vectors ${\bf u}_{k+1}$ and ${\bf v}_{k+1}$ can be found
by using a hyperbolic rotation $H \left( \theta_k \right)$ defined by the
following relations:
\begin{eqnarrayabc}
\label{hypsin}
\sin \theta_k & = & {\bf e}^T_{k+1} {\bf v}_k / {\bf e}^T_k {\bf u}_k\;,
\label{eq35a}\\
\cos \theta_k & = & \sqrt{1- \sin^2 \theta_k}\;,
\label{hypcos}
\end{eqnarrayabc}
and
\be
\left[ \begin{array}{l}
{\bf u}_{k+1}^T\\
{\bf v}_{k+1}^T
\end{array} \right] = H \left( \theta_k \right) \left[ \begin{array}{l}
{\bf u}_k^TZ^T\\
{\bf v}_k^T
\end{array} \right]\;.
\label{newuvhyp}
\ee
The elementary downdating problem has a unique
solution (up to obvious sign changes) if
\[
|{\bf e}^T_k {\bf u}_k| > |{\bf e}^T_{k+1} {\bf v}_k|\;.
\]

The calculation of ${\bf u}_{k+1}$, ${\bf v}_{k+1}$ via~(\ref{newuvhyp}) can
be performed in the obvious manner.  Following common usage,
algorithms which perform downdating in this manner
will be referred to as {\it hyperbolic} downdating algorithms.

Some computational advantages may be obtained
by rewriting~(\ref{elementdown}) as follows:
\[
\left[{\bf u}_{k+1} \; {\bf v}_{k}\right]
\left[ \begin{array}{c}{\bf u}_{k+1}^T\\ {\bf v}_{k}^T\end{array} \right] =
\left[Z{\bf u}_{k} \; {\bf v}_{k+1}\right]
\left[ \begin{array}{c} {\bf u}_{k}^TZ^T\\ {\bf v}_{k+1}^T\end{array} \right]\;.
\]
Consider now an orthogonal rotation matrix $G(\theta_k)$,
\[
G(\theta_k)=\left[ \begin{array}{cc}\cos \theta_k&\sin \theta_k\\
-\sin \theta_k&\cos \theta_k \end{array}\right]\;,
\]
where $\cos \theta_k$ and $\sin \theta_k$ are defined
by~(\ref{hypcos}) and~(\ref{hypsin}), respectively.
Then it is easy to check that
\be
G(\theta_k)\left[ \begin{array}{c}{\bf u}_{k+1}^T\\ {\bf v}_{k}^T\end{array}
\right] =
\left[ \begin{array}{c} {\bf u}_{k}^TZ^T\\ {\bf v}_{k+1}^T\end{array}
\right]\;,
\label{downort1}
\ee
or, equivalently,
\be
\left[ \begin{array}{c}{\bf u}_{k+1}^T\\ {\bf v}_{k}^T\end{array}
\right] =G(\theta_k)^T \left[ \begin{array}{c} {\bf u}_{k}^TZ^T\\ {\bf
v}_{k+1}^T\end{array} \right]\;.
\label{downort2}
\ee
Thus, we may rewrite~(\ref{downort2}) as
\begin{eqnarrayabc}
{\bf v}_{k+1} & = & ({\bf v}_k - \sin \theta_k Z{\bf u}_k)/ \cos
\theta_k\;,
\label{newvort} \\
{\bf u}_{k+1} & = & - \sin \theta_k {\bf v}_{k+1} + \cos \theta_k Z {\bf
u}_k\;.
\label{newuort}
\end{eqnarrayabc}
Note that equation~(\ref{newvort}) is the same as
the second component of~(\ref{newuvhyp}).
However,~(\ref{newuort}) differs from the first component of~(\ref{newuvhyp})
as it uses ${\bf v}_{k+1}$ in place of ${\bf v}_{k}$
to define ${\bf u}_{k+1}$.
It is possible to construct an
alternative algorithm by using the first component of~(\ref{downort1}) to
define ${\bf u}_{k+1}$. This leads to the following formulas:
\begin{eqnarrayabc}
\label{newuort2}
{\bf u}_{k+1} & = & (Z{\bf u}_k-\sin \theta_k {\bf v}_{k})/ \cos
\theta_k\;,\\
{\bf v}_{k+1} & = & - \sin \theta_k {\bf u}_{k+1} + \cos \theta_k {\bf
v}_k\;.
\label{newvort2}
\end{eqnarrayabc}
We call algorithms based on~(\ref{newvort})--(\ref{newuort})
or~(\ref{newuort2})--(\ref{newvort2})
{\it mixed} elementary downdating algorithms.
The reason for considering
mixed algorithms is that they have superior stability
properties to hyperbolic algorithms in the following sense.

Let $\tilde{{\bf u}}_k$, $\tilde{{\bf v}}_k$ be the values of ${\bf
u}_k$, ${\bf v}_k$ that are computed in floating point arithmetic with
relative machine precision $\epsilon$. The computed values
$\tilde{{\bf u}}_k$, $\tilde{{\bf v}}_k$ satisfy a perturbed version
of~(\ref{elementdown}), that is,
\be
\tilde{{\bf u}}_{k+1} \tilde{{\bf u}}^T_{k+1} - \tilde{{\bf v}}_{k+1}
\tilde{{\bf v}}^T_{k+1} = Z \tilde{{\bf u}}_k \tilde{{\bf u}}_k^T Z^T -
\tilde{{\bf v}}_k \tilde{{\bf v}}_k^T + \epsilon G_k + O (\epsilon^2)\;,
\label{elementdownperturb}
\ee
where the second order term $O(\epsilon^2)$ should be understood as a matrix
whose elements are bounded by a constant multiple of $\epsilon^2||G_k||$.
The norm of the perturbation $G_k$ depends on the precise
specification of the algorithm used. It can be shown~\cite{BBDH}
that the term $G_k$ satisfies
\be
\| G_k \| \leq c_m  \left( \| Z{\bf u}_k \|^2+ \|
{\bf v}_k \|^2+ \| {\bf u}_{k+1} \|^2 + \|  {\bf v}_{k+1} \|^2 \right)
\label{normGort}
\ee
when a mixed downdating strategy is used (here $c_m$ is a positive constant).
When hyperbolic downdating is used the term $G_k$ satisfies
\be
\| G_k \| \leq c_h \| H( \theta_k) \| \left( \|
Z{\bf u}_k \| + \| {\bf v}_k \|\right)\left(\| {\bf u}_{k+1} \| + \|
{\bf v}_{k+1} \|\right)\;,
\label{normGhyp}
\ee
where $c_h$ is a positive constant~\cite{BBDH}.
(The constants $c_m$ and $c_h$ are
dependent on implementation details,
but are of order unity and independent of $n$.)
Note the presence of the multiplier~$\| H( \theta_k) \|$ in the
bound~(\ref{normGhyp}) but not in~(\ref{normGort}).
In view of~(\ref{Heigenval}), $\| H( \theta_k) \|$ could be large.
The significance of the multiplier $\| H( \theta_k) \|$ depends on the
context in which the downdating arises. We consider the
implications of the bounds~(\ref{normGort}) and~(\ref{normGhyp})
in Section~\ref{Sec:5} after we make a
connection between downdating and the factorization of Toeplitz matrices.

It is easily seen that a single step of the hyperbolic or mixed downdating
algorithm requires $4(n-k)+O(1)$ multiplications.
A substantial increase in efficiency can be achieved
by considering the following modified downdating problem.  Given
$\alpha_k$, $\beta_k \in \Re$ and ${\bf w}_k$, ${\bf x}_k \in \Re^n$ that
satisfy
\[
{\bf e}^T_j {\bf w}_k =0 \; , \; j <k \;\;\;\;{\rm and}\;\;\;\;
{\bf e}_j^T {\bf x}_k =0 \; , \; j \leq k \;,
\]
find $\alpha_{k+1}$, $\beta_{k+1}$ and ${\bf w}_{k+1}$, ${\bf x}_{k+1}
\in \Re^n$ that satisfy
\[
\alpha^2_{k+1} {\bf w}_{k+1} {\bf w}^T_{k+1} - \beta^2_{k+1} {\bf
x}_{k+1} {\bf x}^T_{k+1} = \alpha^2_k Z{\bf w}_k {\bf w}^T_k Z^T -
\beta^2_k {\bf x}_k {\bf x}^T_k\;,
\]
with
\[
{\bf e}^T_j {\bf w}_k =0 \; , \; j <k \;\;\;\;{\rm and}\;\;\;\;
{\bf e}_j^T {\bf x}_k =0 \; , \; j \leq k \;.
\]
If we make the identification
\begin{eqnarray*}
{\bf u}_k  =  \alpha_k {\bf w}_k\;\;\;\;{\rm and}\;\;\;\;
{\bf v}_k  =  \beta_k {\bf x}_k\;,
\end{eqnarray*}
then we find that the modified elementary downdating problem is equivalent to
the elementary downdating problem.  However, the extra parameters
can be chosen judiciously to eliminate some multiplications.
For example, if we take $\alpha_k = \beta_k$, $\alpha_{k+1} =
\beta_{k+1}$, then
from~(\ref{eq35a}), (\ref{hypcos}) and~(\ref{newuvhyp}),
\begin{eqnarrayabc}
\sin \theta_k & = & {\bf e}^T_{k+1} {\bf x}_k/ {\bf e}^T_k {\bf w}_k\;,
\label{eq320a}\\
\alpha_{k+1} & = & \alpha_k/ \cos \theta_k\;,
\label{eq320b}
\end{eqnarrayabc}
and
\begin{eqnarrayabc}
{\bf w}_{k+1} & = & Z{\bf w}_k - \sin \theta_k {\bf x}_k\;,
\label{scalhypa}\\
{\bf x}_{k+1} & = & - \sin \theta_k Z {\bf w}_k + {\bf x}_k\;.
\label{scalhypb}
\end{eqnarrayabc}
Equations~(\ref{eq320a})--(\ref{scalhypb})
form a basis for a {\it scaled hyperbolic} elementary downdating
algorithm which requires
$2(n-k) + O(1)$ multiplications. This is about half the
number required by the unscaled
algorithm based on~(\ref{newuvhyp}).
(The price is an increased likelihood of underflow or overflow,
but this can be avoided if suitable precautions are taken in the code.)

Similarly, from~(\ref{newvort}) and~(\ref{newuort})
we can obtain a {\it scaled mixed} elementary downdating algorithm via
\begin{eqnarray*}
\sin \theta_k & = & \beta_k {\bf e}^T_{k+1} {\bf x}_k/ \alpha_k {\bf
e}^T_k {\bf w}_k\;,\\
\alpha_{k+1} & = & \alpha_k \cos \theta_k\;,\\
\beta_{k+1} & = & \beta_k/ \cos \theta_k\;,
\end{eqnarray*}
and
\begin{eqnarray*}
{\bf x}_{k+1} & = & {\bf x}_k - \frac{\sin \theta_k \alpha_k}{\beta_k} Z
{\bf w}_k\;,\\
{\bf w}_{k+1} & = & - \frac{\sin \theta_k \beta_{k+1}}{\alpha_{k+1}}
{\bf x}_{k+1} + Z {\bf w}_k\;.
\end{eqnarray*}

The stability properties of scaled mixed algorithms are similar to those
of the corresponding unscaled algorithms~\cite{Fletcher}.

\section{Symmetric Factorization}
\label{Sec:4}
\set

We adopt the following definition from~\cite{Kailath}.

\noindent
{\bf Definition~4.1:} An $n \times n$ symmetric matrix $T$ has
displacement rank 2 iff there exist vectors
${\bf u}$, ${\bf v} \in \Re^n$ such that
\be
T-ZTZ^T= {\bf u}{\bf u}^T-{\bf v}{\bf v}^T \;.
\label{displacement}
\ee
\hspace*{\fill}{\large $\Box$}

The vectors ${\bf u}$ and ${\bf v}$ are called the generators of $T$
and determine the matrix $T$ uniquely.
Whenever we want to stress the dependence of $T$ on ${\bf u}$ and ${\bf v} $ we
write $T=T({\bf u}\:,\: {\bf v})$.

In the sequel we will be concerned with a subset $\cal T$ of all matrices
satisfying~(\ref{displacement}). The subset is defined as follows.

\noindent
{\bf Definition~4.2:} A matrix $T$ is in $\cal T$ iff
\begin{description}
\item [(a)] $T$ is positive definite,
\item [(b)] $T$ satisfies~(\ref{displacement}) with generators ${\bf u}$
and ${\bf v}$,
\item [(c)] ${\bf v}^T {\bf e}_1 =0$, i.e., the first component of ${\bf v}$
is zero.
\end{description}
\hspace*{\fill}{\large $\Box$}

It is well known that positive definite $n \times n$ Toeplitz matrices
form a subset of $\cal T$. Indeed, if $T=(t_{|i-j|})_{i,j=0}^{n-1}$, then
\[
T-ZTZ^T= {\bf u}{\bf u}^T-{\bf v}{\bf v}^T \;,
\]
where
\begin{eqnarray*}
{\bf u}^T & = & \left( t_0\:,\: t_1 \:,\: \ldots \: , \: t_{n-1}
\right)/\sqrt{t_0}\;,\\
{\bf v}^T & = & \left( 0\:,\: t_1 \:,\: \ldots \: , \: t_{n-1}
\right)/\sqrt{t_0}\;.
\end{eqnarray*}
\noindent The set $\cal T$ also contains matrices which are not Toeplitz,
as the following example shows.

\noindent
{\bf Example:} Let
\[
T=\left[ \begin{array}{ccc} 25&20&15\\20&32&29\\15&29&40 \end{array}\right]
\;\;\;,\;\;\;{\bf u}=\left[ \begin{array}{c}5\\4\\3 \end{array}\right]\;\;\;
{\rm and}\;\;\;{\bf v}=\left[ \begin{array}{c}0\\3\\1 \end{array}\right]\;.
\]
It is easy to check that $T$ is positive definite. Moreover,
\[
T-ZTZ^T=\left[ \begin{array}{ccc} 25&20&15\\20&7&9\\15&9&8
\end{array}\right] =
\left[ \begin{array}{ccc}25&20&15\\20&16&12\\15&12&9\end{array}\right]-
\left[ \begin{array}{ccc}0&0&0\\0&9&3\\0&3&1\end{array}\right]=
{\bf u}{\bf u}^T-{\bf v}{\bf v}^T\;.
\]
Hence $T=T({\bf u}\:,\:{\bf v})\in \cal T$, but $T$ is not Toeplitz.

\hspace*{\fill}{\large $\Box$}

We now establish a connection between the elementary downdating
problem and symmetric factorizations of a matrix from the set $\cal T$.

Let $T=T({\bf u}\:,\:{\bf v})\in \cal T$. Set
\[
{\bf u}_1 = {\bf u}, \hspace{.5cm} {\bf v}_1 = {\bf v}
\]
and, for $k=1, \ldots ,n-1,$ solve the elementary downdating problem
defined by~(\ref{elementdown}),
\[
{\bf u}_{k+1} {\bf u}_{k+1}^T - {\bf v}_{k+1} {\bf v}_{k+1}^T = Z {\bf
u_k}{\bf u}^T_kZ^T- {\bf v}_k {\bf v}_k^T\;,
\]
which we assume for the moment has a solution
for each $k$.
On summing over $k = 1,\ldots, n-1$ we obtain
\[
\sum_{k=1}^{n-1}{\bf u}_{k+1}{\bf u}_{k+1}^T-
\sum_{k=1}^{n-1}{\bf v}_{k+1}{\bf v}_{k+1}^T=
\sum_{k=1}^{n-1} Z {\bf u}_k {\bf u}_k^T Z^T-
\sum_{k=1}^{n-1}{\bf v}_k {\bf v}_k^T\;.
\]
If we now observe that, from~(\ref{zerok+1}),
\[
Z {\bf u}_n = {\bf v}_n = 0\;,
\]
we arrive at the following relation:
\be
\sum_{k=1}^{n}{\bf u}_{k}{\bf u}_{k}^T-
Z\left( \sum_{k=1}^{n}{\bf u}_{k}{\bf u}_{k}^T\right) Z^T=
{\bf u}_{1}{\bf u}_{1}^T-{\bf v}_{1}{\bf v}_{1}^T \;,
\label{incalT}
\ee
which implies that $\sum_{k=1}^{n}{\bf u}_{k}{\bf u}_{k}^T\in \cal T$.
Moreover, as matrices having the same generators are identical, we obtain
\[
T = \sum_{k=1}^{n}{\bf u}_{k}{\bf u}_{k}^T = U^TU\;,
\]
where
\[
U = \sum^n_{k=1} {\bf e}_k {\bf u}_k^T
\]
is upper triangular, and hence is the Cholesky factor of $T$.
We have derived, albeit in a rather indirect manner, the basis of an
algorithm for calculating the Cholesky decomposition of a matrix from
the set $\cal T$.

We now return to the question of existence of a solution to the
elementary downdating problem for each $k=1, \ldots , n-1$.  It is easy
to verify that, if $T\in \cal T$, then
$|{\bf e}_1^T {\bf u}_1| > |{\bf e}^T_{2} {\bf v}_1|$.
Using~(\ref{incalT}) and~(\ref{elementdown}),
it can be shown by induction on $k$ that
\[
|{\bf e}_k^T {\bf u}_k| > |{\bf e}^T_{k+1} {\bf v}_k|, \hspace{.5cm} k=2,
\ldots ,n-1.
\]
Consequently, $|\sin \theta_k|<1$ in~(\ref{eq35a}), and
the elementary downdating problem has a solution for
each $k=1, \ldots ,n-1$.

To summarize, we have the following algorithm for factorizing a matrix
$T=T({\bf u}\:,\:{\bf v}) \in \cal T$.

\noindent
{\bf Algorithm} {\it FACTOR(\/$T$)}\/:

\noindent
Set ${\bf u}_1 = {\bf u}$, ${\bf v}_1 = {\bf v}$.\\
For $k=1, \ldots , n-1$  calculate ${\bf u}_{k+1}$, ${\bf
v}_{k+1}$ such that
\begin{eqnarray*}
{\bf u}_{k+1} {\bf u}^T_{k+1} - {\bf v}_{k+1} {\bf v}^T_{k+1} & = & Z
{\bf u}_k {\bf u}^T_k Z^T - {\bf v}_k {\bf v}^T_k\;,\\
{\bf e}^T_{k+1} {\bf v}_{k+1} & = & 0\;.
\end{eqnarray*}
\noindent
Then $T=U^TU,$ where $ U= \sum^n_{k=1} {\bf e}_k {\bf u}_k^T$.

\hspace*{\fill}{\large $\Box$}

In fact we have not one algorithm but
a class of factorization algorithms, where each
algorithm corresponds to a particular way of realizing the elementary
downdating steps. For example, the connection with the scaled
elementary downdating problem is
straightforward. On making the identification
\be
{\bf u}_k = \alpha_k {\bf w}_k\;\;\;\;{\rm and}\;\;\;\;
{\bf v}_k = \beta_k {\bf x}_k\;,
\label{scaling}
\ee
we obtain
\[
T= W^TD^2W\;,
\]
where
\begin{eqnarray*}
W & = & \sum^n_{k=1} {\bf e}_k {\bf w}^T_k\;,\\
D & = & \sum^n_{k=1} \alpha_k {\bf e}_k {\bf e}^T_k\;.
\end{eqnarray*}

It is clear from Section~\ref{Sec:3} that
Algorithm {\it FACTOR(\/$T$)} requires
$2n^2+O(n)$ multiplications when the unscaled version of elementary
downdating is used,
and $n^2+O(n)$ multiplications when the scaled version of
elementary downdating is used.
However, in the sequel we do not dwell on the precise details of
algorithms.  Using~(\ref{scaling}),
we can relate algorithms based on the scaled
elementary downdating problem to those based on the unscaled elementary
downdating problem. Thus, for simplicity, we
consider only the unscaled elementary downdating algorithms.

\section{Analysis of Factorization Algorithms}
\label{Sec:5}
\set

In this section we present a numerical stability
analysis of the factorization of $T\in \cal T$ via
Algorithm {\it FACTOR(\/$T$)}. The result of the
analysis is applied to the case when the matrix $T$ is Toeplitz.

Let $\tilde{{\bf u}}_k$, $\tilde{{\bf v}}_k$ be the values of ${\bf
u}_k$, ${\bf v}_k$ that are computed in floating point arithmetic with
relative machine relative precision $\epsilon$. The computed quantities
$\tilde{{\bf u}}_k$ and $\tilde{{\bf v}}_k$ satisfy the relations
\be
\tilde{{\bf u}}_k = {\bf u}_k + O(\epsilon), \hspace{.5cm} \tilde{{\bf
v}}_k = {\bf v}_k + O(\epsilon),	\label{eq50}
\ee
and the aim of this section is to provide a first order analysis of the
error.
By a first order analysis we mean that
the error can be bounded by a function which has an asymptotic expansion
in powers of $\epsilon$, but we only consider the first term of
this asymptotic expansion.
One should think of $\epsilon \to 0+$ while the problem
remains fixed~\cite{Miller}.
Thus, in this section (except for Corollary~5.1) we omit functions of $n$
from the ``$O$'' terms in relations such as~(\ref{eq50}) and~(\ref{eq51}).

The computed vectors
$\tilde{{\bf u}}_k$, $\tilde{{\bf v}}_k$ satisfy a
perturbed version~(\ref{elementdownperturb}) of~(\ref{elementdown}).
On summing~(\ref{elementdownperturb})
over $k = 1,\ldots,n-1$ we obtain
\[
\tilde{T} -Z \tilde{T} Z^T = \tilde{{\bf u}}_1\tilde{{\bf u}}^T_1 -
\tilde{{\bf v}}_1 \tilde{{\bf v}}_1^T-(Z \tilde{{\bf u}}_n \tilde{{\bf
u}}_n^TZ^T - \tilde{{\bf v}}_n \tilde{{\bf v}}_n^T)+ \epsilon
\sum^{n-1}_{k=1} G_k+ O(\epsilon^2)\;,
\]
where
\begin{eqnarray*}
\tilde{T} & = & \tilde{U}^T \tilde{U}\;, \\
\tilde{U} & = & \sum^n_{k=1} {\bf e}_k \tilde{{\bf u}}_k^T\;.
\end{eqnarray*}
Since
\[
Z \tilde{{\bf u}}_n = O(\epsilon), \hspace{.5cm} \tilde{{\bf v}}_n =
O(\epsilon)\;,
\]
we find that
\be
\tilde{T} - Z \tilde{T} Z^T = \tilde{{\bf u}}_1 \tilde{{\bf u}}^T_1 -
\tilde{{\bf v}}_1 \tilde{{\bf v}}^T_1 + \epsilon \sum^{n-1}_{k=1} G_k +
O(\epsilon^2)\;.			\label{eq51}
\ee
Now define
\be
\tilde{E} = \tilde{T} -T.	\label{eq52}
\ee
Then, using~(\ref{displacement}),~(\ref{eq51}) and~(\ref{eq52}),
\[
\tilde{E}-Z\tilde{E}Z^T= \tilde{{\bf u}}_1 \tilde{{\bf u}}^T_1- {\bf
uu}^T+\tilde{{\bf v}}_1 \tilde{{\bf v}}^T_1 - {\bf vv}^T+ \epsilon
\sum^{n-1}_{k=1} G_k + O( \epsilon^2)\;.
\]
In a similar manner we obtain expressions for
$Z_j\tilde E Z_j^T-Z_{j+1}\tilde E Z_{j+1}^T$, $j=0,\ldots,n-1$.
Summing over $j$ gives
\be
\tilde{E} =
\sum^{n-1}_{j=0} Z_j \Bigl( (\tilde{{\bf u}}_1 \tilde{{\bf u}}_1^T-
{\bf u}_1 {\bf u}_1^T)+(\tilde{{\bf v}}_1 \tilde{{\bf v}}^T_1 - {\bf
v}_1 {\bf v}^T_1)\Bigr) Z^T_j+ \epsilon \sum^{n-1}_{j=0} \sum^{n-1}_{k=1}
Z_jG_kZ^T_j
+O( \epsilon^2)\;.	\label{eq53}
\ee
We see from~(\ref{eq53}) that the error consists of two parts~-- the first
associated with initial errors and the second associated with the fact
that~(\ref{eq51}) contains an inhomogeneous term.
Now
\begin{eqnarray*}
\| \tilde{{\bf u}}_1 \tilde{{\bf u}}_1^T - {\bf uu}^T \| \leq 2 \| {\bf
u} \| \; \| \tilde{{\bf u}}_1 - {\bf u} \| + O( \epsilon^2)\;, \\
\| \tilde{{\bf v}}_1 \tilde{{\bf v}}_1^T - {\bf vv}^T \| \leq 2 \| {\bf
v} \| \; \| \tilde{{\bf v}}_1 - {\bf v} \| +O( \epsilon^2)\;.
\end{eqnarray*}
Furthermore, from~(\ref{displacement}),
\[
Tr(T)-Tr(ZTZ^T)= \| {\bf u} \|^2 - \| {\bf v} \|^2 >0\;,
\]
and hence
\be
\Bigl\| \sum^{n-1}_{j=0} Z_j ( \tilde{{\bf u}}_1 \tilde{{\bf u}}^T_1 - {\bf
uu}^T+ \tilde{{\bf v}}_1 \tilde{{\bf v}}^T_1 - {\bf vv}^T) Z^T_j \Bigr\| \leq
2n \| {\bf u} \| \Bigl( \| \tilde{{\bf u}}_1 - {\bf u} \| + \| \tilde{{\bf v}}_1
-
{\bf v} \| \Bigr) + O(\epsilon^2)\;.		\label{eq54}
\ee
This demonstrates that initial errors do not propagate unduly.  To
investigate the double sum in~(\ref{eq53}) we require a preliminary result.

\vspace{2ex}

\noindent
{\bf Lemma~5.1} ~For $k=1,2,\ldots,\: n-1$ and $j=0,1,2,\ldots,$
\[
\| Z_j {\bf v}_k \| \leq \| Z_{j+1} {\bf u}_k \|\;.
\]

\hfill {\large $\Box$}

\noindent
{\bf Proof} ~Let
\[
T_k = T- \sum^k_{l=1} {\bf u}_l {\bf u}^T_l
=  \sum^n_{l=k+1} {\bf u}_l {\bf u}^T_l\;.
\]
It is easy to verify that
\[
T_k-ZT_kZ^T = Z {\bf u}_k{\bf u}_k^T Z^T - {\bf v}_k {\bf v}_k^T
\]
and, since $T_k$ is positive semi-definite,
\[
Tr\Bigl(Z_jT_kZ^T_j-Z_{j+1}T_kZ^T_{j+1}\Bigr) =
\| Z_{j+1} {\bf u}_k \|^2 - \| Z_j
{\bf v}_k \|^2 \geq 0 \;.
\]

\hfill {\large $\Box$}

We now demonstrate stability when the mixed version of elementary
downdating is used in Algorithm {\it FACTOR(\/$T$)}.
In this case the inhomogeneous term $G_k$ satisfies a shifted version
of~(\ref{normGort}), that is
\be
\| Z_j G_kZ^T_j \| \leq c_m  \Bigl( \| Z_{j+1} {\bf u}_k \|^2+ \| Z_j
{\bf v}_k \|^2+ \| Z_j {\bf u}_{k+1} \|^2 + \| Z_j {\bf v}_{k+1} \|^2 \Bigr)\;,
			\label{shiftnormGort}
\ee
where $c_m$ is a positive constant.         %

\vspace{2ex}

\noindent
{\bf Theorem~5.1} ~Assume that~(\ref{elementdownperturb})
and~(\ref{shiftnormGort}) hold.	%
Then
\[
\| T- \tilde{U}^T \tilde{U} \| \leq 2n \| {\bf u} \|
\Bigl( \| \tilde{{\bf u}}_1 - {\bf u} \| + \| \tilde{{\bf v}}_1 - {\bf v} \|
\Bigr) +4 \epsilon c_m \sum^{n-1}_{j=0} Tr(Z_jTZ^T_j) + O ( \epsilon^2)\;.
\]

\hfill $\Box$

\noindent
{\bf Proof} ~Using Lemma~5.1,
\[
\| Z_jG_kZ^T_j \| \leq 2c_m  \Bigl( \| Z_{j+1} {\bf u}_k \|^2 + \| Z_j {\bf
u}_{k+1} \|^2 \Bigr).
\]
Furthermore, since
\[
Tr(Z_jTZ_j^T)= \sum^n_{k=1} \| Z_j {\bf u}_k\|^2,
\]
it follows that
\be
\Bigl\| \sum^{n-1}_{j=0} \sum^n_{k=1} Z_j G_k Z^T_j \Bigr\| \leq 4 c_m
\sum^{n-1}_{j=0} Tr (Z_jTZ_j^T)\;.
\label{boundZG}
\ee
The result now follows from~(\ref{eq53}), (\ref{eq54}) and~(\ref{boundZG}).

\hfill {\large $\Box$}

For the hyperbolic version of the elementary downdating algorithms
a shifted version of the
weaker bound~(\ref{normGhyp}) on $G_k$ holds (see~\cite{BBDH}), namely
\be
\| Z_j G_k Z^T_j \| \leq c_h \| H( \theta_k) \| ( \|
Z_{j+1} {\bf u}_k \| + \| Z_j {\bf v}_k \|)(\| Z_j {\bf u}_{k+1} \| + \|
Z_j {\bf v}_{k+1} \|)\;.	\label{shiftGhyp} %
\ee
By Lemma~5.1, this simplifies to
\be	%
\|Z_j G_k Z^T_j \| \leq 4 c_h \| H ( \theta_k) \| \; \| Z_{j+1} {\bf u}_k \|
\; \| Z_j {\bf u}_{k+1} \|\;.
\label{boundGhyp}
\ee
The essential difference between~(\ref{normGort}) and~(\ref{normGhyp})
is the occurence of the
multiplier $\| H ( \theta_k) \|$ which can be quite large.  This term
explains numerical difficulties in applications such as
the downdating of a Cholesky decomposition~\cite{BBDH}.
However, because
of the special structure of the matrix $T$,
it is of lesser importance here,
in view of the following result.

\noindent
{\bf Lemma~5.2} ~For $k=1,2,\ldots ,\: n-1$ and $j=0,1, \ldots ,\: n-k$,
\[
\| H( \theta_k) \| \; \| Z_j {\bf u}_{k+1} \| \leq 2 (n-k-j) \| Z_{j+1}
{\bf u}_k \|.
\]

\hfill $\Box$

{\samepage
\noindent{\bf Proof} ~It is easy to verify from~(\ref{newuvhyp}) that
\[
\frac{1 \mp \sin \theta_k}{\cos \theta_k} \bigl( {\bf u}_{k+1} \mp {\bf
v}_{k+1} \bigr) =Z {\bf u}_k \mp {\bf v}_k\;,
\]
and from~(\ref{hyprotH}) that
\[
\| H ( \theta_k ) \| = \frac{1+ | \sin \theta |}{\cos \theta}\;.
\]
Thus,
\begin{eqnarray*}
\| H ( \theta_k) \| \; \| Z_j {\bf u}_{k+1} \| & \leq & \| H (\theta_k) \|
\; \| Z_j {\bf v}_{k+1} \| + \| Z_{j+1} {\bf u}_k \| + \| Z_j {\bf v}_k
\| \\
& \leq & \| H( \theta_k) \| \; \| Z_{j+1} {\bf u}_{k+1} \| +2 \| Z_{j+1}
{\bf u}_k \|\;,
\end{eqnarray*}
where the last inequality was obtained using Lemma~5.1.  Thus
{\samepage \[
\| H ( \theta_k) \| \; \| Z_j {\bf u}_{k+1} \| \leq 2 \sum^{n-k}_{l=j+1}
\| Z_l {\bf u}_k \|\;,
\]
and the result follows.
}
\hfill $\Box$
}

\noindent
{\bf Remark} ~Lemma~5.2 does not hold for the computed quantities unless we
introduce an $O(\epsilon)$ term. However in a first order analysis
we only need it to hold for the exact quantities.

\vspace{2ex}

\noindent
{\bf Theorem~5.2} ~Assume that~(\ref{elementdownperturb})
and~(\ref{shiftGhyp}) hold. Then	%
\[	%
\| T- \tilde{U}^T \tilde{U} \| \leq 2n \| {\bf u} \| \bigl( \|
\tilde{{\bf u}}_1 - {\bf u} \| + \| \tilde{{\bf v}}_1 - {\bf v} \| \bigr)
+ 8 \epsilon c_h \sum^{n-1}_{j=1} (n-j) Tr(Z_jTZ^T_j)
+ O(\epsilon^2)\;.
\]

\hfill $\Box$

\noindent
{\bf Proof} ~Applying Lemma~5.2 to~(\ref{boundGhyp}) gives
\[	%
\| Z_j G_kZ^T_j \| \leq 8 c_h (n-j-1) \| Z_{j+1} {\bf u}_k \|^2\;,
\]
and hence
\begin{eqnarray}
\Bigl\| \sum^{n-1}_{j=0} \sum^{n-1}_{k=1} Z_j G_k Z^T_j \Bigr\|
& \leq & 8 c_h
\sum^{n-1}_{j=1} \sum^{n-1}_{k=1} (n-j) \| Z_j {\bf u}_k \|^2 \nonumber
\\
& \leq & 8 c_h \sum^{n-1}_{j=1} (n-j) Tr (Z_jTZ^T_j)\;.	%
\label{boundZGhyp}
\end{eqnarray}
The result now follows from~(\ref{eq53}),~(\ref{eq54}) and~(\ref{boundZGhyp}).

\hfill $\Box$

Note that, when $T$ is Toeplitz,
\[
Tr(Z_jTZ^T_j) = (n-j)t_0 \;.
\]
Hence, from Theorems 5.1 and 5.2, we obtain our main result on the
stability of the factorization algorithms based on Algorithm {\it
FACTOR(\/$T$)} for a symmetric positive definite Toeplitz matrix:

{\bf Corollary 5.1} ~ The factorization algorithm {\it FACTOR(\/$T$)}
applied to a symmetric positive definite Toeplitz matrix $T$ produces
an upper triangular matrix $\tilde{U}$ such that
\[
T=\tilde{U}^T\tilde{U} + \Delta T\;,
\]
where $\| \Delta T \| = O(\epsilon t_0 n^2)$ when mixed downdating
is used, and $\| \Delta T \| = O(\epsilon t_0 n^3)$ when hyperbolic
downdating is used.

\hfill {\large $\Box$}

\section{The Connection with the Bareiss algorithm}
\label{Sec:6}
\set

\noindent
In his 1969 paper~\cite{Bareiss}, Bareiss proposed an $O(n^2)$ algorithm for
solving Toeplitz linear systems. For a symmetric Toeplitz matrix $T$,
the algorithm, called a {\it symmetric Bareiss algorithm}
in~\cite{Sweet}, can be expressed as follows. Start with a matrix
$A^{(0)}:=T$ and partition it in two ways:
\[
A^{(0)}=\left(\begin{array}{c}U^{(0)}\\T^{(1)}\end{array}\right)\;\;\;,\;\;
A^{(0)}=\left(\begin{array}{c}T^{(-1)}\\L^{(0)}\end{array}\right)\;,
\]
where $U^{(0)}$ is the first row of $T$ and $L^{(0)}$ is the last row of
$T$.
Now, starting from $A^{(0)}$, compute successively
two matrix sequences $\{A^{(i)}\}$ and $\{A^{(-i)}\}$, $i=1,\dots,n-1,$
according to the relations
\be
A^{(i)}= A^{(i-1)}-\alpha_{i-1}Z_iA^{(-i+1)}\;\;\;,\;\;\;
A^{(-i)}=A^{(-i+1)}-\alpha_{-i+1}Z_i^TA^{(i-1)}\;.
\label{UnsymBareiss}
\ee
For $1\leq i \leq n-1,$ partition $A^{(i)}$ and $A^{(-i)}$ as follows:
\[
A^{(i)}=\left(\begin{array}{c}U^{(i)}\\T^{(i+1)}\end{array}\right)\;\;\;,\;\;
A^{(-i)}=\left(\begin{array}{c}T^{(-i-1)}\\L^{(i)}\end{array}\right)\;,
\]
where $U^{(i)}$ denotes the first $i+1$ rows of $A^{(i)}$, and $L^{(i)}$
denotes
the last $i+1$ rows of $A^{(-i)}$. It is shown in~\cite{Bareiss}
{\samepage that
\begin{description}
\item [(a)] $T^{(i+1)}$ and $T^{(-i-1)}$ are Toeplitz,
\item [(b)] for a proper choice
of $\alpha_{i-1}$ and $\alpha_{-i+1}$, the matrices $L^{(i)}$ and
$U^{(i)}$ are lower and upper trapezoidal, respectively,
and
\item [(c)] with the choice of $\alpha_{i-1}$ and $\alpha_{-i+1}$ as in (b),
the Toeplitz matrix $T^{(-i-1)}$ has zero elements in positions $2,\ldots,i+1$
of its first row, while the Toeplitz matrix $T^{(i+1)}$ has zero elements in
positions $n-1,\ldots,n-i$ of its last row.
\end{description}
}
{\samepage
Pictorially,
$$ A^{(i)}=
\left(\begin{array}{c}U^{(i)}\\ \hline \\[-12pt] T^{(i+1)}\end{array}\right)=
{
\left(\begin{array}{cccccccc}
\times &\cdots & \cdots & \cdots & \cdots & \cdots & \cdots & \times \\
0&\times &&&&&& \times \\
\vdots &\ddots &\times &\times &&&& \vdots \\
0&\cdots &0&\times &\cdots &\cdots &\cdots &\times \\
\hline &&&&&&&\\[-12pt]
\times &0&\cdots &0&\times &\cdots &\cdots &\times \\
\vdots & \ddots &&&\ddots & \ddots && \vdots\\
\vdots & & \ddots &&&\ddots &\ddots & \vdots\\
\times &\cdots &\cdots &\times &0 &\cdots &0& \times \\
\end{array}\right)
}
$$
$$ A^{(-i)}=
\left(\begin{array}{c}T^{(-i-1)}\\ \hline\\[-12pt] L^{(i)}\end{array}\right)=
{
\left(\begin{array}{cccccccc}
\times &0&\cdots &0&\times &\cdots &\cdots &\times \\
\vdots & \ddots &&&\ddots &\ddots&& \vdots\\
\vdots & & \ddots &&&\ddots &\ddots& \vdots\\
\times &\cdots &\cdots &\times &0 &\cdots &0& \times \\
\hline &&&&&&&\\[-12pt]
\times &\cdots &\cdots &\times &\times &0&\cdots& 0 \\
\vdots & & &&&\ddots &\ddots & \vdots\\
\vdots &  &&&&&\ddots& 0\\
\times &\times &\cdots &\cdots&\cdots &\cdots &\cdots &\times \\
\end{array}\right)
}
$$
}

After $n-1$ steps, the matrices
$A^{(n-1)}$ and $A^{(-n+1)}$ are lower and upper triangular, respectively.
At step $i$ only rows $i+1,\ldots,n$ of $A^{(i)}$ and rows $1,2,\ldots,n-i$
of $A^{(-i)}$ are modified; the remaining rows stay unchanged.
Moreover, Bareiss~\cite{Bareiss} noticed
that, because of the symmetry of~$T$,
\be
T^{(i+1)}=J_{i+1}T^{(-i-1)}J_n\;\;\;{\rm and}\;\;\;\alpha_{i-1}=\alpha_{-i+1}
\;,
\label{symmetry}
\ee
Here $J_{i+1}$ and $J_n$ are the reversal matrices of dimension
$(i+1)\times (i+1)$ and $n\times n$ respectively.

Now, taking into account~(\ref{symmetry}), it can be seen that
the essential part of a step of the Bareiss algorithm~(\ref{UnsymBareiss})
can be written as follows:	%
\be
\left(\begin{array}{cccc}
t_{i+2}^{(i+1)}&t_{i+3}^{(i+1)}&\ldots&t_n^{(i+1)}\\
0&t_{i+3}^{(-i-1)}&\ldots&t_n^{(-i-1)}
\end{array}\right)=
\left(\begin{array}{cc}1&-\alpha_{i-1}\\
-\alpha_{i-1}&1\end{array}\right)
\left(\begin{array}{cccc}
t_{i+2}^{(i)}&t_{i+3}^{(i)}&\ldots&t_n^{(i)}\\
t_{i+2}^{(-i)}&t_{i+3}^{(-i)}&\ldots&t_{n}^{(-i)}
\end{array}\right)\;,
\label{simple}
\ee
where $(t_{i+2}^{(-i)},t_{i+2}^{(-i)},\ldots ,t_n^{(-i)})$
are the last $(n-i-1)$ components of the first row of $T^{(-i)}$, and
$(t_{i+2}^{(i)},t_{i+3}^{(i)},\ldots ,t_n^{(i)})$ are the last
$(n-i-1)$ components of the first row of
$T^{(i)}$.

Note that~(\ref{simple}) has the same form
as~(\ref{scalhypa})--(\ref{scalhypb}),
and hence a connection between the Bareiss algorithm and algorithm
$FACTOR(T)$ is evident. That such a connection exists was
observed by Sweet~\cite{Sweet}, and later by Delosme and Ipsen~\cite{Ipsen}.
Sweet~\cite{Sweet} related a step of the Bareiss algorithm~(\ref{simple})
to a step of Bennett's downdating procedure~\cite{Bennett}.
Next, he derived the $LU$ factorization of a Toeplitz matrix as a
sequence of Bennett's downdating steps. Finally, he estimated the forward
error in the decomposition using Fletcher and Powell's
methodology~\cite{Fletcher}.
This paper generalizes and presents new derivations of the results
obtained in~\cite{Sweet}.

\section{Numerical examples}
\label{Sec:7}
\set

We adopt from~\cite{Jankowski} the following definitions of forward
and backward stability.

{\bf Definition~7.1:} An algorithm for solving the equation
(\ref{eq11}) is
{\em forward stable } if the computed solution $\tilde{x}$ satisfies
\[
||x-\tilde{x}|| \leq c_1(n) \epsilon {\rm cond}(T)||\tilde{x}||\;,
\]
where ${\rm cond}(T) = ||T||\, ||T^{-1}||$ is the condition number of $T$,
and $c_1(n)$ may grow at most as fast as a polynomial in $n$,
the dimension of the system.

{\bf Definition~7.2:} An algorithm for solving the
equation~(\ref{eq11}) is {\em backward stable} if the computed solution
$\tilde{x}$ satisfies
\[
||T\tilde{x} -b|| \leq c_2(n) \epsilon ||T||\, ||\tilde{x}||\;,
\]
where $c_2(n)$ may grow at most as fast as a polynomial in $n$,
the dimension of the system.

It is known that an algorithm (for solving a system of
linear equations) is backward stable iff there exists a matrix
$\Delta T$ such that
\[
(T+\Delta T)\tilde{x} = b \;\;,\;\; ||\Delta T|| \leq c_3(n) \epsilon ||T||\;,
\]
where $c_3(n)$ may grow at most as fast as a polynomial in $n$.

Note that our definitions do not require the perturbation $\Delta T$
to be Toeplitz, even if the matrix $T$ is Toeplitz.
The case that $\Delta T$ is Toeplitz is discussed
in~\cite{Gohberg,Varah2}.	%
The reader is referred to~\cite{Bunch,Golub,Miller} for
a detailed treatment of roundoff analysis for general matrix algorithms.

It is easy to see that backward stability implies forward stability,
but not vice versa.  This is manifested by the size of the residual
vector.

Cybenko~\cite{Cybenko1} showed that the $L_1$ norm of the inverse of
a $n \times n$ symmetric positive definite Toeplitz matrix $T_n$ is
bounded by
\[
\max \Bigl\{ \frac{1}{\prod_{i=1}^{n-1} \cos^2 \theta_i} \:,\:
\frac{1}{\prod_{i=1}^{n-1} (1+\sin \theta_i)} \Bigr\} \leq
\| T_n^{-1}\|_1 \leq
\prod_{i=1}^{n-1} \frac{1+|\sin \theta_i|}
{1-|\sin \theta_i|}\;,
\]
where $\{ -\sin\theta_i \}_{i=1}^{n-1}$
are quantities called {\it reflection coefficients}.
It is not difficult to pick the reflection coefficients
in such a way that the corresponding Toeplitz matrix $T_n$ satisfies
\[
{\rm cond} (T_n) \approx 1/\epsilon \;.
\]
One possible way of constructing a Toeplitz matrix with given reflection
coefficients
${\{ -\sin \theta_i\}}_{i=1}^{n-1}$ is by tracing the elementary downdating
steps backwards.

An example of a symmetric positive definite Toeplitz
matrix that can be made poorly conditioned by suitable choice of
parameters is the {\em Prolate}
matrix~\cite{Slepian,Varah},	%
defined by
\[            		%
t_k = 	\left\{ \begin{array}{ll}
	2\omega					& \mbox{if $k=0$,}\\
	\frac{\sin (2\pi \omega k)}{\pi k}	& \mbox{otherwise},
	\end{array}
	\right.
\]
where $ 0 \leq \omega \leq \frac{1}{2}$.
For small $\omega$ the eigenvalues of the Prolate matrix cluster
around~0 and~1.

We performed numerical tests in which we solved systems of
Toeplitz linear equations using variants of the Bareiss and Levinson
algorithms, and (for comparison) the standard Cholesky method.
The relative machine precision was
$\epsilon = 2^{-53} \approx 10^{-16}.$
We varied the dimension of the system from~10 to~100,
the condition number of the matrix from~1 to~${\epsilon}^{-1}$,
the signs of reflection coefficients, and the right hand side
so the magnitude of the norm
of the solution vector varied from~1 to~${\epsilon}^{-1}$.
In each test we monitored the errors in the decomposition,
in the solution vector, and the size of the residual vector.

Let $x_B$ and $x_L$ denote the solutions computed by the Bareiss and Levinson
algorithms. Also, let $r_B=Tx_B-b$ and $r_L=Tx_L-b$.
Then for the Bareiss algorithms we always observed that
the scaled residual
\[
s_B\equiv \frac{\|r_B\|}{\epsilon\|x_B\|\|T\|}
\]
was of order unity, as small as would be expected
for a backward stable method.
However, we were not able to find
an example which would demonstrate the superiority of
the Bareiss algorithm based on
mixed downdating over the Bareiss algorithm based on
hyperbolic downdating.
In fact, the Bareiss algorithm based on
hyperbolic downdating often gave slightly smaller errors than
the Bareiss algorithm based on mixed downdating.
In our experiments with Bareiss algorithms, neither
the norm of the error matrix in the
decomposition of $T$ nor the residual error in the solution
seemed to depend in any clear way on $n$, although a quadratic or cubic
dependence would be expected from the worst-case error bounds of
Theorems~5.1--5.2 and Corollary~5.1.

For well conditioned systems the Bareiss and Levinson algorithms behaved
similarly, and gave results comparable to results produced by a general
stable method (the Cholesky method). Differences between the
Bareiss and Levinson algorithms were noticeable only for very
ill-conditioned systems and special right-hand side vectors.

For the Levinson algorithm, when the matrix was very
ill-conditioned and the norm of the
solution vector was of order unity (that is, when the norm of the
solution vector did not reflect the ill-conditioning of the matrix),
we often observed that the scaled residual
\[
s_L\equiv \frac{\|r_L\|}{\epsilon\|x_L\|\|T\|} \,,
\]
was as large as $10^5$.
Varah~\cite{Varah} was the first to observe
this behavior of the Levinson
algorithm on the Prolate matrix.
Higham and Pickering~\cite{Pickering} used a search method
proposed in~\cite{Higham} to
generate Toeplitz matrices for which the Levinson algorithm
yields large residual errors.
However, the search never produced $s_L$ larger than $5\cdot 10^5$.
It plausible that $s_L$ is a slowly increasing function of $n$ and
$1/\epsilon$.

Tables~7.1--7.3 show typical behavior of the Cholesky, Bareiss and Levinson
algorithms for ill-conditioned Toeplitz systems of linear equations when
the norm of the solution vectors is of order unity.
The decomposition
error was measured for the Cholesky and Bareiss algorithms by the quotient
$||T-L\cdot L^T||/(\epsilon\cdot ||T||)$, where $L$ was
the computed factor of $T$.
The solution error was measured by the quotient
$||x_{comp}-x||/||x||$, where $x_{comp}$ was the computed solution vector.
Finally, the residual error was measured by the
quotient \mbox{$||T\cdot x_{comp}-b||/(||T||\cdot ||x_{comp}||\cdot \epsilon)$.}

\vspace{2ex}
\pagebreak[3]	%

{\samepage
\begin{center}
\begin{tabular}{l|c|c|c|c}
& decomp. error & soln. error & resid. error \\ \hline \hline
&&&&\\
Cholesky&$5.09\cdot 10^{-1}$&$7.67\cdot 10^{-3}$& $1.25\cdot 10^0$\\
Bareiss(hyp)&$3.45\cdot 10^{0}$&$1.40\cdot 10^{-2}$&$8.72\cdot 10^{-1}$\\
Bareiss(mixed)&$2.73\cdot 10^{0}$&$1.41\cdot 10^0$&$1.09\cdot 10^0$\\
Levinson&&$5.30\cdot 10^0$&$4.57\cdot 10^3$
\end{tabular}\\[1ex]
\nopagebreak
Table 7.1: Prolate matrix, $n=21$, $\omega =0.25$, $cond=3.19\cdot 10^{14}$ \\
\end{center}

\vspace{1ex}
\pagebreak[3]	%

\begin{center}
\begin{tabular}{l|c|c|c|c}
& decomp. error & soln. error & resid. error \\ \hline \hline
&&&&\\
Cholesky&$1.72\cdot 10^{-1}$&$6.84\cdot 10^{-2}$& $3.11\cdot 10^{-1}$\\
Bareiss(hyp)&$2.91\cdot 10^{0}$&$2.19\cdot 10^{-1}$&$1.15\cdot 10^{-1}$\\
Bareiss(mixed)&$3.63\cdot 10^{0}$&$2.48\cdot 10^{-1}$&$2.47\cdot 10^{-1}$\\
Levinson&&$5.27\cdot 10^{-1}$&$1.47\cdot 10^5$
\end{tabular}\\[1ex]\nopagebreak
Table 7.2: Reflection coefficients $|\sin \theta_i|$ of the same magnitude
$|K|$ but\\
alternating signs,
$|K|=0.8956680108101296$, $n=41$, $cond=8.5\cdot10^{15}$
\end{center}

\vspace{1ex}
\pagebreak[3]

\begin{center}
\begin{tabular}{l|c|c|c|c}
& decomp. error & soln. error & resid. error \\ \hline \hline
&&&&\\
Cholesky&$8.51\cdot 10^{-1}$&$3.21\cdot 10^{-2}$& $4.28\cdot 10^{-1}$\\
Bareiss(hyp)&$8.06\cdot 10^{0}$&$1.13\cdot 10^{-1}$&$2.28\cdot 10^{-1}$\\
Bareiss(mixed)&$6.71\cdot 10^{0}$&$1.16\cdot 10^{-1}$&$3.20\cdot
10^{-1}$\\
Levinson&& $2.60\cdot 10^{-1}$&$1.06\cdot 10^5$\\
\end{tabular}\\[1ex]
\nopagebreak
Table 7.3: Reflection coefficients $|\sin \theta_i|$
of the same magnitude but\\
alternating signs,
$|K|=0.9795872473975045$, $n=92$, $cond=2.77\cdot 10^{15}$\\
\end{center}
}
\pagebreak[3]

\section{Conclusions}
\label{Sec:8}
\set

This paper generalizes and presents new derivations of results
obtained earlier by Sweet~\cite{Sweet}. The bound in
Corollary~5.1 for the case of mixed downdating is stronger than
that given in~\cite{Sweet}. %
The applicability of the Bareiss algorithms based
on elementary downdating steps is extended to a class of
matrices, satisfying Definition~4.2, which includes
symmetric positive definite Toeplitz matrices. The approach via
elementary downdating greatly simplifies roundoff error analysis.
Lemmas~5.1 and~5.2 appear to be new. The
stability of the Bareiss algorithms follows directly from
these Lemmas and the results on the roundoff error analysis for
elementary downdating steps given in~\cite{BBDH}.

The approach via downdating can be extended
to the symmetric factorization of positive definite
matrices of displacement rank~$k\ge 2$ (satisfying additional
conditions similar to those listed in Definition~4.2); see~\cite{Kailath}.
For matrices of displacement rank $k$ the factorization algorithm uses
elementary rank-$k$ downdating via hyperbolic Householder or mixed
Householder reflections~\cite{BojStein,Rader}.

We conclude by noting that the Bariess algorithms guarantee small
residual errors in the
sense of Definition~7.2, but the Levinson algorithm
can yield residuals at least
five orders of magnitude larger than those expected for a backward stable
method.
This result suggests that, if the Levinson algorithm is used in
applications where the reflection coefficients are not known in advance
to be positive, the residuals
should be computed to see if they are acceptably small. This can be done
in $O(n \log n)$ arithmetic operations (using the FFT).

It is an interesting open question whether the Levinson algorithm can give
scaled residual errors which are arbitrarily large
(for matrices which are numerically nonsingular).
A related question is whether the Levinson algorithm,
for positive definite Toeplitz matrices
$T$ without a restriction on the reflection coefficients,
is stable in the sense of Definitions~7.1 or~7.2.

\pagebreak[3]


\begin{thebibliography}{99}

\bibitem{Alexander} S.T. Alexander, C-T. Pan and R.J. Plemmons,
``Analysis of a Recursive least Squares Hyperbolic Rotation Algorithm for
Signal Processing'', {\em Linear Algebra and Its Applications}\/,
vol 98, pp 3-40, 1988.

\bibitem{Bareiss} E.H. Bareiss, ``Numerical Solution of Linear Equations
with Toeplitz and Vector Toeplitz Matrices'', {\em Numerische Mathematik}\/,
vol 13, pp 404-424, 1969.

\bibitem{Bennett} J.M. Bennett, ``Triangular factorization of modified
matrices'', {\em Numerische Mathematik,} vol 7, pp 217-221, 1965.

\bibitem{BBH} A.W. Bojanczyk, R.P. Brent and F.R. de Hoog, ``QR
Factorization of Toeplitz Matrices'', {\em Numerische Mathematik}\/, vol 49,
pp 81-94, 1986.

\bibitem{BBH2} A.W. Bojanczyk, R.P. Brent and F.R. de Hoog, ``Stability
Analysis of Fast Toeplitz Linear System Solvers'',
Report CMA-MR17-91, Centre for Mathematical Analysis,
The Australian National University, August 1991.

\bibitem{BBDH} A.W. Bojanczyk, R.P. Brent, P. Van Dooren and F.R. de Hoog,
``A Note on Downdating the Cholesky Factorization'',
{\em SIAM J.\ Sci.\ Statist.\ Comput.}\/,
vol 8, pp 210-220, 1987.

\bibitem{Bojanczyk} A.W. Bojanczyk and A. Steinhardt, ``Matrix Downdating
Techniques for Signal Processing'', {\em Proceedings of the SPIE
Conference on Advanced Algorithms and Architectures for Signal
Processing III}, vol 975, pp 68-75, 1988.

\bibitem{BojStein} A.W. Bojanczyk and A.O. Steinhardt, ``Stabilized Hyperbolic
Householder Transformations'',
{\em IEEE Trans.\ Acoustics, Speech and Signal Processing}\/,
vol ASSP-37, 1989, pp 1286-1288.

\bibitem{Bunch} J.R. Bunch, ``The Weak and Strong Stability of Algorithms in
Numerical Linear Algebra'', {\em Linear Algebra and Its Applications},
vol 88/89, pp 49-66, 1987.


\bibitem{Cybenko1} G. Cybenko, ``The Numerical Stability of the
Levinson-Durbin Algorithm for Toeplitz Systems of Equations'',
{\em SIAM J.\ Sci.\ Statist.\ Comput.}\/,
vol 1, pp 303-319, 1980.


\bibitem{Ipsen} J-M. Delosme and I.C.F. Ipsen, ``From Bareiss's algorithm to
the stable computation of partial correlations'', {\it Journal of
Computational and Applied Mathematics,} vol 27, pp 53-91, 1989.

\bibitem{Fletcher} R. Fletcher and M.J.D. Powell, ``On the Modification
of $LDL^T$ Factorizations'', {\em Mathematics of Computation}\/, vol 28,
pp 1067-87, 1974.

\bibitem{Gohberg} I. Gohberg, I. Koltracht and D. Xiao, ``On the solution of
the Yule-Walker equation'', {\em Proceedings of the SPIE Conference on
Advanced Algorithms and Architectures for Signal
Processing IV}, vol 1566, July 1991.            	%

\bibitem{Golub} G.H. Golub and C. Van Loan, {\em{Matrix Computations}},
second edition, Johns Hopkins Press, Baltimore, Maryland, 1989.

\bibitem{Higham} N.J. Higham, ``Optimization by Direct Search in Matrix
Computations'',
Numerical Analysis Report No.\ 197, University of Manchester, England,
1991; to appear in {\it SIAM J. Matrix Anal. Appl.}

\bibitem{Pickering} N.J. Higham and R.L. Pickering, private communication.

\bibitem{Jankowski} M. Jankowski and H. Wozniakowski, ``Iterative Refinement
Implies Numerical Stability'', {\em BIT}\/, vol 17, pp 303-311, 1977.

\bibitem{Kailath} T. Kailath, S.Y. Kung and M. Morf, ``Displacement Ranks of
Matrices and Linear Equations'', {\em J. Math. Anal. Appl.}\/, vol 68, pp
395-407, 1979.


\bibitem{Miller} W. Miller and C. Wrathall,
{\em Software for Roundoff Analysis of Matrix Algorithms},
Academic Press, 1980.

\bibitem{Rader} C.M. Rader and A.O. Steinhardt, ``Hyperbolic Householder
Transformations'', {\em IEEE Transaction on Acoustics, Speech and Signal
Processing,} vol ASSP-34, 1986, pp 1584-1602.

\bibitem{Slepian} D. Slepian, ``Prolate spheroidal wave functions, Fourier
analysis, and uncertainty V: the discrete case'',
{\em Bell Systems Tech. J.}\/, vol 57, 1978, pp 1371-1430.

\bibitem{Sweet} D. Sweet, ``Numerical Methods for Toeplitz Matrices'', PhD
thesis, University of Adelaide, 1982.

\bibitem{Varah} J.M. Varah, ``The Prolate Matrix: A Good Toeplitz Test
Example'', {\em SIAM Conference on Control, Signals and Systems},
San Francisco, 1990.
Also {\em Linear Algebra Appl.}, vol 187, 1993, pp 269-278.

\bibitem{Varah2} J.M. Varah, ``Backward Error Estimates for Toeplitz and
Vandermonde Systems'', preprint, 1992.
Also Tech.\ Report 91-20, Univ.\ of British Columbia, Sept.\ 1991.

\end{thebibliography}
\end{document}